\documentclass[10pt,twoside]{amsart}
\usepackage{lmodern}
\usepackage[T1]{fontenc}
\usepackage{amsmath,amsfonts,amsthm,amssymb,amscd}
\usepackage{mathrsfs}
\usepackage{euscript}
\usepackage{latexsym}
\usepackage[pdftex]{graphicx}
\usepackage{enumitem}
\usepackage{bbm}
\usepackage{accents}

\usepackage{geometry}
\newtheorem{lemma}{Lemma}
\newtheorem{theorem}{Theorem}

\usepackage{hyperref}  
\hypersetup{colorlinks=true, linktoc=all,
    linkcolor=blue,
		citecolor=blue,
		urlcolor=cyan,
		bookmarksopen=true
		}

\title{Short quadratic exponential  sums with primes in major arcs}
\author{Zarullo Rakhmonov, Rahbaroi Hotami, Mashraf Sharifzoda}
\address{A.Dzhuraev Institute of Mathematics,  National Academy of Sciences of Tajikistan}
\email{zarullo.rakhmonov@gmail.com, zarullo-r@rambler.ru}
\date{}
\begin{document}

\begin{abstract}
Using the second moment of Dirichlet $L$-functions on the critical line over the major arcs ${\mathfrak M}({\mathscr L}^b)$, with
$\tau=y^3x^{-1}{\mathscr L}^{-b_1}$, and excluding a small neighborhood of the centers of these arcs, i.e., those $\alpha$ satisfying $|\alpha-\frac{a}{q}|>(8\pi y^2)^{-1}$, for $y\ge x^{1-\frac{\scriptstyle1}{\scriptstyle{9-4\sqrt{2}}}}{\mathscr L}^{c_2}$, where
$c_2=\frac{2A+24+(\sqrt{2}-1)b_1}{2\sqrt{2}-1}$, we obtain the estimate
$$
S_2(\alpha;x,y)=\sum_{x-y<n\le x}\Lambda(n)e(\alpha n^2)\ll y\,{\mathscr L}^{-A}.
$$
Moreover, in a small neighborhood of the center of the major arcs, defined by $|\alpha-\frac{a}{q}|\le(8\pi y^2)^{-1}$, for $y\ge x^{\frac{\scriptstyle5}{\scriptstyle8}}{\mathscr L}^{1,5A+0,25b+18}$, an asymptotic formula with a remainder term is obtained for $S_2(\alpha;x,y)$, where $A$, $b_1$, and $b$ are arbitrary fixed positive constants, and ${\mathscr L}=\ln xq$.
\end{abstract}

\maketitle

\section{Introduction}
According to Dirichlet's theorem on the approximation of real numbers by rational numbers, every $\alpha$ in the interval $[-1/\tau, 1-1/\tau]$,  can be represented in the form
\begin{equation*}
 \alpha =\frac{a}{q}+\lambda ,\quad (a,q)=1,\qquad 1\le q\le \tau , \qquad |\lambda |\le \frac{1}{q\tau }.
\end{equation*}
By $\mathfrak{M}(P)$ we denote those numbers $\alpha$ for which $q\le P$, and by $\mathfrak{m}(P)$ we denote the remaining $\alpha$. $\mathfrak{M}(P)$ and $\mathfrak{m}(P)$ are called the major and minor arcs, respectively.

Here is the English translation of the LaTeX mathematical text:

J.~Liu and  T.~Zhan~\cite{Liu+Zhan-ActaArith-1996,Liu+Zhan-Mh.Math1999}, studying the sum $S_2(\alpha;x,y)$, $y\ge x^{\theta_2}$, conditionally (under the assumption of the validity of the extended Riemann hypothesis) for $\theta_2=\frac23+\varepsilon$ and unconditionally, using M. Jutila's theorem~\cite{Jitila-1991} on the fourth moment of Dirichlet $L$-functions on the critical line for $\theta_2=\frac{11}{16}+\varepsilon$, proved that for any $A>0$ there exist constants $c_i>0$, $i=1,2,3$, and the relation
\begin{align*}
S_2(\alpha;x,y)=\sum_{x-y<n\le x}\hspace{-10pt}\Lambda(n)e(\alpha n^2)=\left\{
                  \begin{array}{ll}
                    M_2(\alpha;x,y)+O\left(\dfrac{y}{{\mathscr L}^A}\right),\quad & q\le{\mathscr L}^{c_1},\ |\lambda|\le\dfrac1{xy{\mathscr L}^{c_2}}; \\
                    O\left(\dfrac{y}{{\mathscr L}^A}\right), & \hbox{ otherwise,}
                  \end{array}                \right.\\
M_2(\alpha;x,y)=\frac1{\varphi(q)}\sum_{\substack{h=1\\ (h,q)=1}}^qe\left(-\frac{ah^2}{q}\right)\int_{x-y}^xe(\lambda u^2)du,\qquad \tau=\frac{y^3}{x{\mathscr L}^{c_3}}.
\end{align*}
holds. Then G.~S.~Liu and H.~X.~Lao~\cite{Lu+Lao-2007-Mon.fur Math} unconditionally proved that $\theta_2=\frac23+\varepsilon$.

Using the theorem on the second moment of Dirichlet $L$-functions on the critical line~\cite{Zhan Tao-Acta-Math-Sinica-1992}, or more precisely its corollary on the density of nontrivial zeros of Dirichlet $L$-functions in narrow rectangles of the critical strip, the first author~\cite{RakhmonovZKh-Arxiv-20261112,RakhmonovZKh-Cheb-Sb-2021-22-4} proved for short trigonometric sums $S_k(\alpha;x,y)$ for $k \ge 3$ an asymptotic formula with a remainder term in a small neighborhood of the center of the major arcs and obtained a nontrivial estimate on the major arcs themselves, except in a small neighborhood of their centers (see also~\cite{RakhmonovZKh+RFZ-TrMIAN-296-2016,RakhmonovZKh+RFZ-DANRT-59-2016, RakhmonovZKh+SAA-FPM-DANRT-63(5-6)-2020, RakhmonovZKh+SAA-FPM-DANRT-63(7-8)-2020,HRL+ShMS-DANRT-2021-64-1-2}).

The methodology for studying the behavior of the sums $S_k(\alpha; x, y)$ for $k = 2$ differs in some aspects from the case $k \ge 3$. In this paper, taking into account the peculiarities of the quadratic sum and using the method presented in~\cite{RakhmonovZKh-Arxiv-20261112}, for short quadratic trigonometric sums with primes $S_2(\alpha; x, y)$ in a small neighborhood of the center of the major arcs ${\mathfrak M}({\mathscr L}^b)$ for $\tau = y^3 x^{-1} {\mathscr L}^{-b_1}$, an asymptotic formula with a remainder term is obtained. On the major arcs ${\mathfrak M}({\mathscr L}^b)$, except in a small neighborhood of their centers, a nontrivial estimate is proved.
{\theorem\label{Teorema povedenie kor trigsumm s prost chislami} Let $x\ge x_0$, $A$, $b_1$, $b$ be arbitrary fixed positive numbers, $1\le q\le{\mathscr L}^b$,
$$
\alpha=\frac{a}{q}+\lambda,\quad (a,q)=1,\quad |\lambda|\le\frac{1}{q\tau},\quad \tau=\frac{y^3}{x{\mathscr L}^{b_1}}.
$$
Then for $|\lambda|\le\left(8\pi y^2\right)^{-1}$ and $y\ge x^{\frac{\scriptstyle5}{\scriptstyle8}}{\mathscr L}^{1.5A+0.25b+18}$ the equality
$$
S_2(\alpha;x,y)=\frac{1}{\varphi(q)}\sum_{\substack{n=1\\(a,q)=1}}^qe\left(\frac{an^2}{q}\right)\int_{x-y}^xe(\lambda u^2)du+O\left(y{\mathscr L}^{-A}\right),
$$
holds, and for $|\lambda|>\left(8\pi y^2\right)^{-1}$ and $y\ge x^{1-\frac{\scriptstyle1}{\scriptstyle{9-4\sqrt{2}}}}{\mathscr L}^{c_2}$, $c_2=\frac{2A+24+\left(\sqrt{2}-1\right)b_1}{2\sqrt{2}-1}$, the estimate
$$
S_2(\alpha;x,y)\ll y{\mathscr L}^{-A}.
$$ holds.}

This theorem is weaker than in the above-mentioned papers~\cite{Liu+Zhan-Mh.Math1999,Lu+Lao-2007-Mon.fur Math}, and its proof uses the theorem on the second moment of Dirichlet $L$-functions on the critical line~\cite{Zhan Tao-Acta-Math-Sinica-1992}.

\section{Proof of Theorem \ref{Teorema povedenie kor trigsumm s prost chislami}}
\subsection{Reduction of the proof of the theorem to the estimate of the sum $W(\alpha;x,y)$}
Proceeding analogously to the proof of the main theorem of~\cite{RakhmonovZKh-Cheb-Sb-2021-22-4}, setting $k = 2$, for $x \ge x_0$ and arbitrary fixed positive numbers $A$, $b$, $b_1$, for every $\alpha$ from the interval $[-\varrho, 1-\varrho]$, where $\varrho\tau = 1$, we have the representation:
$$
\alpha = \frac{a}{q} + \lambda,\quad (a, q) = 1,\quad 1 \le q \le \tau,\quad |\lambda| \le \frac{1}{q\tau},
$$
where $1 \le q \le {\mathscr L}^b$ and $\tau = y^3 x^{-1} {\mathscr L}^{-b_1}$. Then for
$$
T_0 = \left(x y^{-1} + |\lambda| x^2\right) q^{\frac{1}{2}} {\mathscr L}^{A + 3}
$$
we have the asymptotic formula:
\begin{align}
S_2(\alpha;x,y)&=\frac{1}{\varphi(q)}\sum_{\substack{n=1\\(a, q)=1}}^qe\left(\frac{a n^2}{q}\right)\int_{x - y}^xe(\lambda u^2)\,du-W(\alpha;x,y)+O\left(\frac{y}{{\mathscr L}^A}\right), \label{formula Sk(alpha;x,y)=int+W1+W+R2-3} \\
W(\alpha;x,y)&\ll\frac{{\mathscr L}}{\sqrt{q}}\sum_{\chi\bmod q}{\EuScript{W}}(\lambda,\chi),\qquad{\EuScript{W}}(\lambda,\chi) =\sum_{\substack{|\gamma|\le T_0\\\rho\ne\beta_1}}\left|I(\rho,\lambda)\right|, \label{Formula W(lambda,x,y)<<..sum W(lambda,chi)-3} \\
|I(\rho,\lambda)|&\ll x^\beta\min\left(\frac{y}{x},\, \frac{1}{\min|\gamma+4\pi\lambda u^2|},\,\frac1{\sqrt{|4\pi\lambda u^2-\gamma|}}\right), \label{Formula otsenka I(rho,lambda)-3}
\end{align}
where $\beta_1$ is a real zero, if modulo $q$ there exists a real character $\chi_1$ such that $L(s, \chi_1)$ has a real zero $\beta_1 \ge 1 - c / \ln q$.

Also, in the main theorem of~\cite{RakhmonovZKh-Cheb-Sb-2021-22-4} it was shown that in formulas~(\ref{Formula W(lambda,x,y)<<..sum W(lambda,chi)-3}) and~(\ref{Formula otsenka I(rho,lambda)-3}) we may assume without loss of generality that $\lambda \ge 0$. The case $\lambda \le 0$ reduces to the case $\lambda \ge 0$ by means of the relation ${\EuScript{W}}(\lambda, \chi) = {\EuScript{W}}(\chi, -\lambda)$.

\subsection{Derivation of the asymptotic formula for $S_2(\alpha;x,y)$ in the case $\lambda\le\left(8\pi y^2\right)^{-1}$}
This case is proved in exactly the same way as in item <<3.2. Proof of the theorem in the case $\lambda\le\left(2\pi k^2x^{k-2}y^2\right)^{-1}$>> of~\cite{RakhmonovZKh-Cheb-Sb-2021-22-4}, if we set $k=2$.

\subsection{Estimate of $W(\alpha;x,y)$ in the case $\lambda>\left(8\pi y^2\right)^{-1}$}
We introduce an additional parameter $H=8\pi \lambda xy$, which for $\lambda>\left(8\pi y^2\right)^{-1}$ satisfies the condition
$$
H\ge\frac xy.
$$
Everywhere below, without loss of generality, we assume that the parameter $y$ satisfies the condition
\begin{equation}\label{Formula neravenstva y}
y=x^{1-\frac1{9-4\sqrt{2}}}{\mathscr L}_x^{c_2}.
\end{equation}
We partition all nontrivial zeros $\rho=\beta+i\gamma$ of the function $L(s,\chi)$ with the condition $|\gamma|\le T_0 $ into the sets $D_1$, $D_2$ and $D_3$:
\begin{align*}
&D_1=\left\{\rho:\ -T_0\le\gamma<-4\pi \lambda x^2-H\right\},  \\
&D_2=\left\{\rho:\ -4\pi \lambda x^2-H\le\gamma\le-4\pi \lambda(x-y)^2+H\right\}, \\
& D_3=\left\{\rho:\ -4\pi \lambda(x-y)^2+H<\gamma\le T_0 \right\}.
\end{align*}
By ${\EuScript{W}}_j(\lambda,\chi)$, $j=1,\ 2,\ 3$ we denote the sum of the moduli of the integral $I(\rho,\lambda)$ over the zeros $\rho$ belonging to the set $D_j$. Using this notation, we represent the sum ${\EuScript{W}}(\lambda,\chi)$, which is defined by formula (\ref{Formula W(lambda,x,y)<<..sum W(lambda,chi)-3}), in the form
\begin{align}\label{Formula W(chi,lambda)=W1+W2+W3}
&{\EuScript{W}}(\lambda,\chi)={\EuScript{W}}_1(\lambda,\chi)+{\EuScript{W}}_2(\lambda,\chi)+{\EuScript{W}}_3(\chi,\lambda).
\end{align}

\subsubsection{Estimate of ${\EuScript{W}}_1(\lambda,\chi)$}  Adding the term $4\pi \lambda u^2$, $x-y\leq u \le x$, to all three terms of the inequality by means of which the set $D_1$ is defined, we obtain
\begin{align*}
&D_1=\left\{\rho:\ -T_0+4\pi \lambda u^2\le\gamma+6\pi \lambda u^2<-4\pi \lambda x^2+4\pi \lambda u^2-H\right\},  \\
\end{align*}
The function $4\pi \lambda u^2$ is monotonically increasing on the interval $x-y\leq u\le x$, therefore for the right boundary of the set $D_1$ we have
$$
-4\pi \lambda x^2+4\pi \lambda u^2-H\le-H.
$$
Therefore, if $\rho\in D_1$, then the inequality $\gamma+4\pi \lambda u^2<-H$ holds. Consequently, for the monotonically increasing function $\gamma+4\pi \lambda u^2$ on the interval $x-y\le u\le x$ the relation
\begin{align*}
&\min|\gamma+4\pi \lambda u^2|=-\max(\gamma+4\pi \lambda u^2)=-\gamma-4\pi \lambda x^2\ge H,
\quad \text{if}  \quad \rho\in D_1,
\end{align*}
holds. From this, using the second estimate of formula (\ref{Formula otsenka I(rho,lambda)-3}), we find
$$
{\EuScript{W}}_1(\lambda,\chi)\ll\sum_{\rho \in D_1}\frac{x^\beta }{-\gamma-4\pi \lambda x^2}.
$$
We partition all nontrivial zeros in the set
$$
D_1=\left\{\rho:\ -T_0\le\gamma<-4\pi \lambda x^2-H\right\}=\left\{\rho:\ H\le-\gamma-4\pi \lambda x^2<T_0-4\pi \lambda x^2\right\},
$$
into classes  $D_{11},\ \ldots,\ D_{1t}$,  $t\ll T_0H^{-1}$. We assign to the class $D_{1n}$ those zeros $\rho$ for which the condition
\begin{align*}
nH<&-\gamma-4\pi \lambda x^2\le(n+1)H,\hspace{26pt}\text{ if}\quad 1\le n\le t,
\end{align*}
holds. Consequently we obtain
\begin{align*}
{\EuScript{W}}_1(\lambda,\chi)&\ll\sum_{n=1}^t\sum_{\rho\in D_{1n}}\frac{x^\beta}{-\gamma-4\pi \lambda x^2}\le\sum_{n=1}^t\sum_{\rho\in D_{1n}}\frac{x^\beta }{nH}\le\\
&\le\frac{{\mathscr L}}H\max_{1\le n\le t}\sum_{\rho\in D_{1n}}x^{\beta}\le\frac{{\mathscr L}}H\max_{|T|\le T_0}\sum_{T-H<\gamma\le T}x^\beta.
\end{align*}
Using the relation $H=8\pi\lambda xy$, we represent this estimate in the form
\begin{align*}
{\EuScript{W}}_1(\lambda,\chi)\ll\frac{{\mathscr L}}{\lambda xy}\max_{|T|\le T_0}\sum_{T-\lambda x^2y<\gamma\le T}x^\beta.
\end{align*}

\subsubsection{Estimate of ${\EuScript{W}}_3(\lambda,\chi)$}
Adding the term $4\pi \lambda u^2$, $x-y\leq u\le x$, to all three terms of the inequality defining the set $D_3$, we obtain
\begin{align*}
& D_3=\left\{\rho:\ 4\pi \lambda u^2-4\pi \lambda(x-y)^2+H<\gamma+4\pi \lambda u^2\le T_0+4\pi \lambda u^2 \right\}.
\end{align*}
The function $4\pi\lambda u^2$ is monotonically increasing on the interval $x-y\leq u\le x$, therefore for the left boundary of the set $D_3$ we have
$$
4\pi \lambda u^2-4\pi \lambda(x-y)^2+H\ge H.
$$
Thus, if $\rho\in D_3$, then the inequality $\gamma+4\pi \lambda u^2>H$ holds. Consequently, for the monotonically increasing function $\gamma+4\pi \lambda u^2$ on the interval $x-y\le u\le x$ the relation
\begin{align*}
&\min|\gamma+4\pi \lambda u^2|=\min(\gamma+4\pi \lambda u^2)=\gamma+4\pi \lambda(x-y)^2\ge H,
\quad \text{if}\quad\rho\in D_3,
\end{align*}
holds. From this and from the second estimate (\ref{Formula otsenka I(rho,lambda)-3}), we find
$$
{\EuScript{W}}_3(\lambda,\chi)\ll\sum_{\rho \in D_3}\frac{x^\beta }{\gamma+4\pi \lambda(x-y)^2}.
$$
We partition all nontrivial zeros in the set
\begin{align*}
 D_3&=\left\{\rho:\ H\le\gamma+4\pi \lambda(x-y)^2<T_0+4\pi \lambda(x-y)^2\right\},
\end{align*}
into classes  $D_{31},\ \ldots,\ D_{3t},$, $t\ll T_0H^{-1}$. We assign to the class $D_{3n}$ those zeros $\rho$ for which the condition
\begin{align*}
nH<\gamma+4\pi \lambda(x-y)^2\le(n+1)H,\hspace{26pt}\text{ if}\quad 1\le n\le t,
\end{align*}
holds. Therefore, proceeding analogously to the case of the estimate of the sum ${\EuScript{W}}_1(\lambda,\chi)$, we have
\begin{align*}
{\EuScript{W}}_3(\lambda,\chi)\ll&\sum_{n=1}^t\sum_{\rho\in D_{3n}}\frac{x^\beta}{\gamma+4\pi \lambda(x-y)^2}\le\sum_{n=1}^t\sum_{\rho\in D_{3n}}\frac{x^\beta }{nH}\le\frac{{\mathscr L}}H\max_{1\le n\le t}\sum_{\rho\in D_{3n}}x^{\beta}\le\\
&\le\frac{{\mathscr L}}H\max_{|T|\le T_0}\sum_{T-H<\gamma\le T}x^\beta\ll\frac{{\mathscr L}}{\lambda xy}\max_{|T|\le T_0}\sum_{T-\lambda x^2y<\gamma\le T}x^\beta.
\end{align*}

\subsubsection{Estimate of ${\EuScript{W}}_2(\lambda,\chi)$}
Using the explicit value of the parameter $H=8\pi \lambda xy$, we introduce the following notation
\begin{equation}\label{formula opr T1 T2,H1}
\begin{split}
&T_1=4\pi\lambda x^2+H=4\pi\lambda x^2+8\pi\lambda xy,\\
&H_1=4\pi \lambda x^2-4\pi \lambda(x-y)^2+2H=24\pi\lambda xy-4\pi\lambda y^2,\\
&T_1-H_1=4\pi\lambda x^2-16\pi\lambda xy+4\pi\lambda y^2=\pi\lambda x^2-16\pi\lambda xy\left(1+O\left(\frac yx\right)\right),\\
&T_2=T_1+4\pi\lambda u^2=4\pi\lambda(x^2+u^2+2xy),
\end{split}
\end{equation}
and represent the set $D_2$ in the form
\begin{equation}\label{Formula mnozhestva D2}
\begin{split}
D_2&=\left\{\rho:\ -4\pi \lambda x^2-H\le\gamma\le-4\pi \lambda(x-y)^2+H\right\}=\\
&=\left\{\rho:\ T_1-H_1\le-\gamma\le T_1\right\}=\\
&=\left\{\rho:\ T_2-H_1\le4\pi \lambda u^2-\gamma\le T_2\right\}.
\end{split}
\end{equation}
From this, using the conditions $x-y<u\le x$ and \eqref{Formula neravenstva y}, we obtain
\begin{align}
T_2-H_1=&4\pi\lambda(x^2+u^2-4xy+y^2)\ge 4\pi\lambda(x^2+(x-y)^2-4xy+y^2)=\nonumber\\
&=8\pi\lambda(x^2-3xy+y^2)\gg\lambda x^2.\label{formula otsenka T2-H1 snizu}
\end{align}
Thus, if $\rho$ belongs to the set $D_2$, then using the third estimate of formula (\ref{Formula otsenka I(rho,lambda)-3}), the representation of the set $D_2$ in the form (\ref{Formula mnozhestva D2}) through $T_2$, the estimate (\ref{formula otsenka T2-H1 snizu}), and then again the representation of the set $D_2$ in the form (\ref{Formula mnozhestva D2}) but this time through $T_1$, we successively obtain
\begin{align*}
|{\EuScript{W}}_2&(\lambda,\chi)=\sum_{\rho\in D_2}|I(\rho,\lambda)|=\sum_{T_2-H_1\le4\pi\lambda u^2-\gamma\le T_2}|I(\rho,\lambda)|\ll\\
&\ll\sum_{T_2-H_1\le4\pi\lambda u^2-\gamma\le T_2}\frac{x^\beta}{\sqrt{4\pi\lambda u^2-\gamma}}\le\sum_{T_2-H_1\le4\pi\lambda u^2-\gamma\le T_2}\frac{x^\beta}{\sqrt{T_2-H_1}}\le\\
&\ll\sum_{T_2-H_1\le4\pi\lambda u^2-\gamma\le T_2}\frac{x^\beta}{\sqrt{\lambda x^2}}=\sum_{T_1-H_1\le-\gamma\le T_1}\frac{x^\beta}{\sqrt{\lambda x^2}}=\sum_{T_1-H_1\le\gamma\le T_1}\frac{x^\beta}{\sqrt{\lambda x^2}}.
\end{align*}
The boundaries of the last sum over $\gamma$, that is, the quantities $T_1$ and $T_1-H_1$, according to their definitions in formula (\ref{formula opr T1 T2,H1}) are quantities of order $\lambda x^2$. The length of this sum, that is, the parameter $H_1$, is also, according to formula (\ref{formula opr T1 T2,H1}), a quantity of order $24\pi\lambda xy\asymp\lambda xy$. Splitting in this sum the summation interval $T_1-H_1\le\gamma\le T_1$ into a finite number of intervals of the form
$$
\kappa\lambda x^2-\lambda x^2y<\gamma\le \kappa\lambda x^2,
$$
where the constant $\kappa$ takes a value from the interval
$$
4\pi \left(1+\frac{2y}x-\frac{y^2}{x^2}\right)<\kappa\le4\pi \left(1+\frac{4y}x\right),
$$
we obtain
$$
{\EuScript{W}}_2(\lambda,\chi)\ll\frac{1}{\sqrt{\lambda x^2}}\sum_{\kappa\lambda x^2-\lambda xy<\gamma\le\kappa\lambda x^2}x^\beta.
$$
In this estimate, without loss of generality for convenience we assume that $\kappa=1$, and also using the inequality $\lambda x^2<T_0$, we have
$$
{\EuScript{W}}_2(\lambda,\chi)\ll\frac{1}{\sqrt{\lambda x^2}}\max_{|T|\le T_0}\sum_{T-\lambda xy<\gamma\le T}x^\beta.
$$

\subsubsection{Reduction of the estimate of $W(\alpha;x,y)$ to the estimates ${\mathscr B}_1$ and ${\mathscr B}_2$}
Substituting this estimate and the estimates ${\EuScript{W}}_1(\lambda,\chi)$ and ${\EuScript{W}}_2(\lambda,\chi)$ into (\ref{Formula W(chi,lambda)=W1+W2+W3}), and then using the condition $\lambda>\left(8\pi y^2\right)^{-1}$, we obtain:
\begin{align*}
{\EuScript{W}}(\lambda,\chi)&\ll\left(\frac1{\sqrt{\lambda x^2}}+\frac{{\mathscr L}}{\lambda xy}\right)\max_{|T|\le T_0}\sum_{T-\lambda xy<\gamma\le T}x^\beta \ll\frac{{\mathscr L}}{\sqrt{\lambda qx^2}}\max_{|T|\le T_0}\sum_{T-\lambda xy<\gamma\le T}x^\beta.
\end{align*}
From this and from the estimate (\ref{Formula W(lambda,x,y)<<..sum W(lambda,chi)-3}), we find
$$
|W(\alpha;x,y)|\ll\frac{{\mathscr L}^2}{\sqrt{\lambda qx^2}}\max_{|T|\le T_0}\sum_{\chi\bmod q}\sum_{T-\lambda xy<\gamma\le T}x^\beta= \frac{{\mathscr L}^2}{\sqrt{\lambda qx^2}}\max_{|T|\le T_0}{\EuScript{V}}_q(T,\lambda xy).
$$
Using for ${\EuScript{V}}_q(T,U)$ with $U=\lambda xy$ formula (11) of~\cite{RakhmonovZKh-Cheb-Sb-2021-22-4}, we obtain
\begin{align}\label{Formula otsenka summi W'}
|W(\alpha;x,y)|&\ll\frac{{\mathscr L}^2}{\sqrt{\lambda qx^2}}\max_{\substack{|T|\le T_0\\ 0.5\le u\le1-\delta}}x^u\sum_{\chi\bmod q}\left(N(u,T,\chi)-N(u,T-\lambda xy,\chi)\right),\\
&\delta=\delta(q,T_0)=\frac{c_1}{\max\left(\ln q,\ \left(\ln(T_0+3)\ln\ln(T_0+3)\right)^\frac{3}{4}\right)}. \nonumber
\end{align}
From the explicit value of the parameter $T_0$, the condition $\lambda>\left(8\pi y^2\right)^{-1}$ and the relation (\ref{Formula neravenstva y}), we find
\begin{align*}
\frac{T}{(\lambda xy)^3}\le&\frac{T_0}{(\lambda xy)^3}=\frac{\left(xy^{-1}+|\lambda|x^2\right)q^{\frac12}{\mathscr L}^{A+3}}{(\lambda xy)^3}= \frac{\sqrt{q}{\mathscr L}^{A+3}}{\lambda^3x^2y^4}+ \frac{\sqrt{q}{\mathscr L}^{A+3}}{\lambda^2xy^3}\le\\
 &\le\frac{(8\pi)^2y^2{\mathscr L}^{A+0.5b+3}}{x^2}+\frac{(8\pi)^2y{\mathscr L}^{A+0.5b+3}}{x}\ll{\mathscr L}^{-7}.
\end{align*}
From this it follows that in the sum over $\chi\hspace{-5pt}\mod q$ in (\ref{Formula otsenka summi W'}) the condition $\lambda xy\ge T^\frac13$ holds, that is, the following lemma can be applied to this sum.

{\lemma\label{Lemma plotnostn-teor Zhan Tao}{\rm \cite{Zhan Tao-Acta-Math-Sinica-1992}}. Let $\varepsilon$ be an arbitrarily small positive constant, ` and $T^{\frac{35}{108}+\varepsilon}\le H\le T$, then the estimates
\begin{align*}
&\sum_{\chi\hspace{-7pt}\mod q}\left(N(u,T+H,\chi)-N(u,T,\chi)\right)\ll\left\{
                                                       \begin{array}{ll}
                                                         (qH)^{\frac{4}{3-2u}(1-u)}(\ln qH)^9, & \hbox{for} \quad\dfrac{1}{2}\le u\le\dfrac{3}{4} \vspace{7pt} \\
                                                         (qH)^{\frac{2}{u}(1-u)+\varepsilon}, & \hbox{for}\quad\dfrac{3}{4}\le u\le 1,
                                                       \end{array}
                                                     \right.
 \end{align*}}
hold.

Setting in lemma \ref{Lemma plotnostn-teor Zhan Tao}
\begin{equation}\label{formula epsilon vib B2}
\varepsilon=\min\left(\frac{53}{3591\cdot 10^4},\ \frac{\delta}3\right),
\end{equation}
we have
\begin{align}
&|W(\alpha;x,y)|\ll{\mathscr B}_1+{\mathscr B}_2,\label{Formula otsenka summi B1}\\
&{\mathscr B}_1=\frac{{\mathscr L}^{12}}{\sqrt{\lambda qx^2}}\max_{0.5\le u\le0.75}x^u\left(q\lambda xy\right)^{\frac{4-4u}{3-2u}},\nonumber \\
&{\mathscr B}_2=\frac{{\mathscr L}^2}{\sqrt{q\lambda x^2}}\max_{0.75\le u\le1-\delta(q,T_0)}x^u\left(q\lambda xy\right)^{\frac{2}{u}(1-u)+\varepsilon}. \nonumber
\end{align}

\subsubsection{Estimate of ${\mathscr B}_1$}
We have
\begin{align*}
&{\mathscr B}_1= \frac{y^\frac12{\mathscr L}^{12}}{x^\frac12}\max_{0.5\le u\le0.75}\exp\left(f_1(u)\right),\\
&f_1(u)=u\ln x+\left(\frac{1}{u-1.5}+\frac32\right)\ln\left(q\lambda xy\right), \\
&f_1'(u)=\ln x-\frac{\ln\left(q\lambda xy\right)}{(u-1.5)^2},\qquad f_1''(u)=\frac{2\ln\left(q\lambda xy\right)}{(u-1.5)^3}<0.\hspace{87pt}
\end{align*}
Let us show that the inequality $f_1'(0.5)>0$ holds. Using the conditions
$$
0\le\lambda\le\frac1{q\tau},\qquad \tau=\frac{y^3}{x{\mathscr L}^{b_1}}, \qquad y\ge x^{1-\frac1{9-4\sqrt{2}}}{\mathscr L}^{c_2},
$$
we find
\begin{align*}
-f_1'(0.5)&=\ln\left(q\lambda y\right)\le\ln\left(\frac{y}{\tau}\right)=\ln\left(\frac{x{\mathscr L}^{b_1}}{y^2}\right)\le  \ln\left(\frac{{\mathscr L}^{b_1}}{x^{1-\frac2{9-4\sqrt{2}}}{\mathscr L}^{2c_2}}\right)= \\
&=\ln\left(x^{-\frac{17}{31+8\sqrt{2}}}{\mathscr L}^{b_1-2c_2}\right)=-\frac{17}{31+8\sqrt{2}}{\mathscr L}+(b_1-4c_2)\ln{\mathscr L}<0.
\end{align*}

\textbf{Case $q\lambda xy\le x^\frac9{16}$}. Using the condition of the case under consideration, we have
$$
f_1'(u)\ge f_1'(0.75)=\ln x-\frac{16}{9}\ln(q\lambda xy)\ge\ln x-\frac{16}{9}\ln x^\frac9{16}=0,
$$
consequently the function $f_1(u)$ is increasing on the interval $0.5\le u\le0.75$. Using this property, and then the relation (\ref{Formula neravenstva y}), we find
\begin{align*}
&{\mathscr B}_1=\frac{y^\frac12{\mathscr L}^{12}}{x^\frac12}\exp\left(f_1(0.75)\right)=\frac{y^\frac12{\mathscr L}^{12}}{x^\frac12}\exp\left(\frac{3}{4}\ln x+\frac16\ln(q\lambda xy)\right)=\\
&=\left(q\lambda xy\right)^\frac16x^\frac14y^\frac12{\mathscr L}^{12}<x^\frac{11}{32}y^\frac12{\mathscr L}^{12}= y\left(\frac{x^\frac{11}{16}{\mathscr L}^{2A+24}}{y}\right)^\frac12{\mathscr L}^{-A}\ll y{\mathscr L}^{-A}.
\end{align*}

\textbf{Case $q\lambda xy>x^\frac9{16}$}. Introducing additional parameters $\mu$ and $c$, denoting the length of the sum $y$ in the form
\begin{equation}\label{Formula znachenie y}
  y=x^{1-\mu^2}{\mathscr L}^c,
\end{equation}
we seek the largest value of $\mu$ and the smallest value of $c$ for which the estimate
\begin{equation}\label{Formula otsenka summi W'1-2sl}
{\mathscr B}_1\ll  y{\mathscr L}^{-A}
\end{equation}
holds. Using the conditions of the case under consideration, namely the relations
$$
\frac1{yx^\frac7{16}}<q\lambda\le\frac1{\tau},\qquad \tau=\frac{y^3}{x{\mathscr L}^{b_1}},\qquad y\ge x^{1-\frac1{9-4\sqrt{2}}}{\mathscr L}^{c_2}>x^\frac12{\mathscr L}^{-0.5b_1},
$$
we find
$$
q\lambda xy\le\frac{xy}{\tau}=\left(\frac{x}{y}\right)^2{\mathscr L}^{b_1}< \left({\mathscr L}^{-0.5b_1}\right)^2{\mathscr L}^{b_1}=1.
$$
Therefore
\begin{align*}
&f_1'(0.5)=\ln x-\ln(q\lambda xy)>\ln x>0,
&f_1'(0.75)=\ln x-\frac{16}{9}\ln(q\lambda xy)<0,
\end{align*}
and the point
$$
u_0=\frac32-\left(\frac{\ln(q\lambda xy)}{\ln x}\right)^\frac12<\frac32-\left(\frac9{16}\right)^\frac12=\frac34,
$$
where $f_1'(u_0)=0$, belongs to the interval $[0.5,\ 0.75]$. On this interval, as we have already noted, the inequality $f''(u)<0$ also holds. Therefore on this interval $0.5\le u\le 0.75$ the graph of the function $f(u)$ is convex upward, consequently
\begin{align}
&{\mathscr B}_1=\frac{y^\frac12{\mathscr L}^{12}}{x^\frac12}\exp\left(f_1(u_0)\right)=
\frac{y^\frac12{\mathscr L}^{12}}{x^\frac12}\exp\left(u_0\ln x+\left(\frac32+\frac{1}{u_0-1.5}\right)\ln\left(q\lambda xy\right)\right)\nonumber\\
&=\frac{y^\frac12{\mathscr L}^{12}}{x^\frac12}\exp\left(\left(\frac{3}{2}-\left(\frac{\ln(q\lambda xy)}{\ln x}\right)^\frac12\right)\ln x\hspace{-3pt}+\hspace{-4pt}\left(\frac32-\left(\frac{\ln x}{\ln(q\lambda xy)}\right)^\frac12\right)\ln\left(q\lambda xy\right)\right)\nonumber\\
&=xy^\frac12{\mathscr L}^{12}\exp\left(\frac32\ln\left(q\lambda xy\right)-2\left(\ln(q\lambda xy)\ln x\right)^\frac12\right) =xy^\frac12{\mathscr L}^{12}\exp\left(g\left(q\lambda xy\right)\right),\label{Formula otsenka summi B11}
\end{align}
where $g(t)=\frac32\ln t-2\left(\ln t\ln x\right)^\frac12$. From the conditions
$$
q\lambda xy\ge x^\frac9{16},\qquad \lambda\le\frac1{q\tau},\qquad \tau=\frac{y^3}{x{\mathscr L}^{b_1}},
$$
and formula (\ref{Formula znachenie y}) it follows that
\begin{align}\label{Formula interval izm qlambdax..}
&x^\frac9{16}\le q\lambda xy\le\frac{xy}{\tau}=\frac{x^2{\mathscr L}^{b_1}}{y^2}=x^u{\mathscr L}^v, \qquad
u=2\mu^2, \quad v=b_1-2c.
\end{align}
Therefore, taking into account the relation $t=q\lambda xy\ge x^\frac9{16}$, we have
\begin{align*}
&g'(t)=\frac3{2t}-\frac{(\ln x)^\frac12}{t(\ln t)^\frac12}
=\frac{9\ln t-4\ln x}{2t(\ln t)^\frac12(3(\ln t)^\frac12+2(\ln x)^\frac12)}>0,
\end{align*}
that is, $g(t)$ is an increasing function on its interval of variation, therefore
$$
g\left(q\lambda x^2y\right)\le g\left(x^u{\mathscr L}^v\right).
$$
Using this inequality and the relation
$$
y^{-\frac12}=x^{-\frac12+\frac{u}4}{\mathscr L}^{-\frac{b_1-v}4},
$$
which is a consequence of relation (\ref{Formula interval izm qlambdax..}), we estimate the right-hand side of (\ref{Formula otsenka summi B11}) as follows
\begin{align}
{\mathscr B}_1&\le xy^\frac12{\mathscr L}^{12}\exp\left(g\left(x^u{\mathscr L}^v\right)\right)=
yx^{\frac12+\frac{u}{4}}{\mathscr L}^{12-\frac{b_1-v}{4}}\exp\left(g\left(x^u{\mathscr L}^v\right)\right)=\nonumber\\
&=yx^{\frac12+\frac{u}{4}}{\mathscr L}^{12-\frac{b_1-v}{4}}\exp\left(\frac32\ln\left(x^u{\mathscr L}^v\right)- 2\left(\ln\left(x^u{\mathscr L}^v\right){\mathscr L}\right)^\frac12\right)=\nonumber\\
&=y{\mathscr L}^{12-\frac{b_1-v}4+\frac{3v}2}\cdot x^{\frac12+\frac{u}4+\frac{3u}2}\exp\left(-2\left(u{\mathscr L}^2+v{\mathscr L}\ln{\mathscr L}\right)^\frac12\right)=\nonumber\\
&=y{\mathscr L}^{12-\frac{b_1}{4}+\frac{7v}{4}}x^{\frac{7u}{4}+\frac12} \exp\left(-2\sqrt{u}{\mathscr L}\left(1+\frac{v\ln{\mathscr L}}{u{\mathscr L}}\right)^\frac12\right).\label{Formula otsenka summi W'1}
\end{align}
Using for $|t|\le 0.1$ the formula
$$
(1+t)^\frac12=1+\frac{t}{2}+R_1(t),\qquad R_1(t,\theta)=-\frac{(1-\theta)t^2}{8(1+\theta t)^{\frac32}},
$$
which is a consequence of the expansion of the function $f(t)=(1+t)^\frac12$ in the Maclaurin series with a remainder term in Cauchy form
\begin{align*}
f(t)=1+&\sum_{n=1}^m\frac{\frac12\left(\frac12-1\right)\ldots\left(\frac12-(n-1)\right)}{n!}t^n+R_m(t,\theta),\\ &R_m(t,\theta)=\frac{f^{(m+1)}(\theta t)}{(m+1)!}(1-\theta)^mt^{m+1},\qquad 0<\theta<1,\\
&f'(t)=\frac{1}{2}(1+t)^{-\frac12},\qquad f''(t)=-\frac{1}{4}(1+t)^{-\frac32},
\end{align*}
with $m=1$. In this formula, considering the remainder term $R_1(t,\theta)$ as a function of $\theta$, we have
\begin{align*}
&\frac{\partial R_1(t,\theta)}{\partial\theta}=\frac{t^2(1+\theta t)^\frac32+(1-\theta)t^2\cdot\frac32(1+\theta t)^\frac12t}{8(1+\theta t)^2}
=\frac{((3-\theta)t+2)t^2}{16(1+\theta t)^\frac52}>0.
\end{align*}
Therefore
$$
\min R_1(t,\theta)=R_1(t,0)=-\frac{t^2}{8}, \qquad (1+t)^\frac12=1+\frac{t}{2}+R_1(t)\ge 1+\frac{t}{2}-\frac{t^2}{8}.
$$
Consequently
\begin{align*}
&-2\sqrt{u}{\mathscr L}\left(1+\frac{v\ln{\mathscr L}}{u{\mathscr L}}\right)^\frac12\le -2\sqrt{u}{\mathscr L}\left(1+\frac{v\ln{\mathscr L}}{2u{\mathscr L}}-\frac{v^2\ln^2{\mathscr L}}{8u^2{\mathscr L}^2}\right)=\\
&=\ln x^{-2\sqrt{u}}+\ln{\mathscr L}^{-\frac{v}{\sqrt{u}}}+\frac{v^2\ln^2{\mathscr L}}{4u^\frac32{\mathscr L}}\le \ln x^{-2\sqrt{u}}+\ln{\mathscr L}^{-\frac{v}{\sqrt{u}}}+1.
\end{align*}
Substituting the found estimate into the right-hand side of (\ref{Formula otsenka summi W'1}), we find
\begin{align}
{\mathscr B}_1\ll&y{\mathscr L}^{12-\frac{b_1}{4}+\frac{7v}{4}}x^{\frac{7u}{4}+\frac12}\cdot x^{-2\sqrt{u}} {\mathscr L}^{-\frac{v}{\sqrt{u}}}=yx^{\frac74u-2\sqrt{u}+\frac12} {\mathscr L}^{12-\frac{b_1}4+\left(\frac74-\frac{1}{\sqrt{u}}\right)v}.\label{Formula otsenka summi W'1-2}
\end{align}
Using the relations in formula (\ref{Formula interval izm qlambdax..}), we express the exponents of $x$ and ${\mathscr L}$ in the right-hand side of the last inequality in terms of the parameters $\mu$ and $c$, denoting them respectively by $\kappa(\mu)$ and $\omega(\mu,c)$; we have
\begin{align*}
\kappa(\mu)&=\frac72\mu^2-2\sqrt{2}\mu+\frac12=\frac12\left(7\mu^2-4\sqrt{2}\mu+1\right),\\
\omega(\mu,c)&=12-\frac{b_1}4+\left(\frac74-\frac1{\sqrt{2}\mu}\right)(b_1-2c)=\\
&=12+\left(\frac32-\frac1{\sqrt{2}\mu}\right)b_1-\left(\frac72-\frac{\sqrt{2}}{\mu}\right)c.
\end{align*}
The quantity $\kappa=\kappa(\mu)$ is a quadratic polynomial in $\mu$ having two positive roots. The larger root is $\mu_2$:
\begin{align}
\mu_2=\frac{2\sqrt{2}+1}7=\frac1{2\sqrt{2}-1},\qquad \mu_2^2=&\frac1{9-4\sqrt{2}}. \label{Formula znachenie mu_k^2}
\end{align}
Therefore the largest value of the parameter $\mu$ for which $\kappa=\kappa(\mu)$---the exponent of $x$ in the estimate (\ref{Formula otsenka summi W'1-2})---is equal to zero, is the root $\mu_2$. Moreover $\omega(\mu,c)$---the exponent of ${\mathscr L}$---for $\mu=\mu_2$ is determined as follows
\begin{align*}
\omega(\mu_2,c)&=12+\left(\frac32-\frac{2\sqrt{2}-1}{\sqrt{2}}\right)b_1-\left(\frac72-\sqrt{2}(2\sqrt{2}-1)\right)c=\\
&=12+\frac{\sqrt{2}-1}2b_1-\frac{2\sqrt{2}-1}2c.
\end{align*}
and for $c=c_2$, where
\begin{align*}
&c_2=\frac{2\left(A+12+\frac{\sqrt{2}-1}2b_1\right)}{2\sqrt{2}-1}=
 \frac{2A+24+\left(\sqrt{2}-1\right)b_1}{2\sqrt{2}-1},
\end{align*}
the equality $\omega(\mu_k,c_2)=-A$ holds, that is, the estimate (\ref{Formula otsenka summi W'1-2}) turns into the estimate (\ref{Formula otsenka summi W'1-2sl}). Thus, from the representations of the parameter $y$ in the form (\ref{Formula znachenie y}) and of the parameter $\mu_k^2$ in the form (\ref{Formula znachenie mu_k^2}) it follows that the estimate (\ref{Formula otsenka summi W'1-2sl}) holds for
$$
y\ge x^{1-\frac1{9-4\sqrt{2}}}{\mathscr L}^{c_2}.
$$

\subsubsection{Estimate of ${\mathscr B}_2$}
We represent the sum ${\mathscr B}_2$ in the form
\begin{align*}
{\mathscr B}_2=&\frac{{\mathscr L}^2}{\sqrt{q\lambda x^2}}\max_{0.75\le u\le1-\delta(q,T_0)}x^u\left(q\lambda xy\right)^{\frac{2}{u}(1-u)+\varepsilon}= \frac{y^\frac12{\mathscr L}^2}{x^\frac12}\left(q\lambda xy\right)^{-\frac12}\max_{0.75\le u\le1-\delta}f_2(u), \\
&f_2(u)=x^u\left(q\lambda xy\right)^{\frac{2}{u}(1-u)+\varepsilon},\qquad \delta=\delta(q,T_0),
\end{align*}
where the function $f_2(u)$ and its second-order derivative $f''(u)$ are positive:
\begin{align*}
f_2'(u)=&f_2(u)\left(\ln x-\frac2{u^2}\ln\left(q\lambda xy\right)\right),\\
f_2''(u)=&f_2(u)\left(\left(\ln x-\frac2{u^2}\ln\left(q\lambda xy\right)\right)^2+\frac4{u^2}\ln\left(q\lambda xy\right)\right)\ge\\ &\ge\frac{4 f_2(u)}{u^2}\ln\left(q\lambda xy\right)\ge\frac{4 f_2(u)}{u^2}\ln\left(\frac{qx}{8\pi y}\right)>0.
\end{align*}
Consequently the graph of the function $f_2(u)$ is convex downward, therefore
\begin{align*}
{\mathscr B}_2\le&\frac{y^\frac12{\mathscr L}^2}{x^\frac12}\left(q\lambda xy\right)^{-\frac12}\left(f_2(0.75)+f_2(1-\delta)\right)=\\
&=\frac{y^\frac12}{x^\frac12}\left(x^\frac34\left(q\lambda xy\right)^{\frac16+\varepsilon}+ x^{1-\delta }\left(q\lambda xy\right)^{-\frac12+\frac{2\delta}{1-\delta}+\varepsilon}\right){\mathscr L}^2=\\
&=y\left(x^\frac14y^{-\frac12}\left(q\lambda xy\right)^\frac{\scriptstyle{1+6\varepsilon}}{\scriptstyle6}+y^{-\frac12}x^{\frac12-\delta }\left(q\lambda xy\right)^{-\frac12+\frac{2\delta}{1-\delta}+\varepsilon}\right){\mathscr L}^2.
\end{align*}
Using the condition of the case under consideration for the parameter $\lambda$, that is, the relation
$$
\dfrac1{8\pi y^2}<\lambda\le\dfrac{1}{q\tau}=\dfrac{x{\mathscr L}^{b_1}}{qy^3},
$$
we find
\begin{equation}\label{Formula otsenka summi B2-2}
\begin{split}
{\mathscr B}_2&\ll y\left(x^\frac14y^{-\frac12}\left(\frac{x^2{\mathscr L}^{b_1}}{y^2}\right)^\frac{\scriptstyle{1+6\varepsilon}}{\scriptstyle6}+ y^{-\frac12}x^{\frac12-\delta}\left(\frac{x}y\right)^{-\frac{\scriptstyle1}{\scriptstyle2}+ \frac{\scriptstyle{2\delta}}{\scriptstyle{1-\delta}}+\varepsilon}\right){\mathscr L}^2=\\
&=y\left(\frac{x^{2+\frac{\scriptstyle3}{\scriptstyle{2(1+6\varepsilon)}}}{\mathscr L}^{b_1}}{y^{2+ \frac{\scriptstyle3}{\scriptstyle{1+6\varepsilon}}}}\right)^\frac{\scriptstyle{1+6\varepsilon}}{\scriptstyle6}\hspace{-22pt}{\mathscr L}^2+ y\left(\frac{x^\frac{\scriptstyle{\delta+\delta^2+(1-\delta)\varepsilon}}{\scriptstyle{2\delta+(1-\delta)\varepsilon}}}{y} \right)^{\frac{\scriptstyle{2\delta}}{\scriptstyle{1-\delta}}+\varepsilon}\hspace{-30pt}{\mathscr L}^2=\\
&=y\left(\frac{x^{1-\frac{\scriptstyle3}{\scriptstyle{10+24\varepsilon}}} {\mathscr L}^{\frac{\scriptstyle{b_1+6\varepsilon}}{\scriptstyle{5+12\varepsilon}}}}{y}\right)^\frac{\scriptstyle{5+12\varepsilon}}{\scriptstyle6} \hspace{-22pt}{\mathscr L}^2+\left(\frac{x^\frac{\scriptstyle{\delta+\delta^2+(1-\delta)\varepsilon}}{\scriptstyle{2\delta+(1-\delta)\varepsilon}}}{y} \right)^{\frac{\scriptstyle{2\delta}}{\scriptstyle{1-\delta}}+\varepsilon}\hspace{-30pt}{\mathscr L}^2=\\
&=y\left(\frac{x^{1-\frac{\scriptstyle1}{\scriptstyle{9-4\sqrt{2}}}}{\mathscr L}^{c_2}}y \right)^\frac{\scriptstyle{5+12\varepsilon}}{\scriptstyle6}\hspace{-22pt}x^{g(\varepsilon)}{\mathscr L}^{h(\varepsilon)} +y\left(\frac{x^{1-\frac{\scriptstyle1}{\scriptstyle{9-4\sqrt{2}}}}{\mathscr L}^{c_2}}y \right)^{\frac{\scriptstyle{2\delta}}{\scriptstyle{1-\delta}}+\varepsilon}\hspace{-30pt}x^{f(\varepsilon)}{\mathscr L}^{v(\varepsilon)},
\end{split}
\end{equation}
where the functions $g(\varepsilon)$,  $h(\varepsilon)$, $f(\varepsilon)$ and $v(\varepsilon)$ are determined by the following relations:
\begin{align*}
&g(\varepsilon)=\left(\frac{1}{9-4\sqrt{2}}-\frac{3}{10+24\varepsilon}\right)\frac{5+12\varepsilon}6;\\
&h(\varepsilon)=\left(\frac{b_1+6\varepsilon}{5+12\varepsilon}-c_2\right)\frac{5+12\varepsilon}6+2b;\\
&f(\varepsilon)=\left(\frac{\delta+\delta^2+(1-\delta)\varepsilon}{2\delta+(1-\delta)\varepsilon}-\left(1-\frac1{9-4\sqrt{2}}\right)\right) \left(\frac{2\delta}{1-\delta}+\varepsilon\right);\\
&v(\varepsilon)=-c_2\left(\frac{2\delta}{1-\delta}+\varepsilon\right)+2.
\end{align*}
Now, successively using the inequalities $40\sqrt{2}<57$ and $c_2>0$, as well as the definition of the parameter $\varepsilon$ from (\ref{formula epsilon vib B2}), we estimate each of the functions $g(\varepsilon)$,  $h(\varepsilon)$, $f(\varepsilon)$ and $v(\varepsilon)$ from above. We have
\begin{align*}
g(\varepsilon)
&=-\frac{57-40\sqrt{2}}{588}+\left(3-\frac{57-40\sqrt{2}}{245}\right)\varepsilon<\\
&<-\frac{57-40\sqrt{2}}{588}+3\varepsilon=-\frac1{84(57+40\sqrt{2})}+3\varepsilon< -\frac1{84\cdot 114}+3\varepsilon=\\
&=-\frac1{9576}+3\varepsilon=-\frac1{10000}+3\left(\varepsilon-\frac{53}{3591\cdot 10^4}\right)\le-\frac1{10000}.
\end{align*}
\begin{align*}
h(\varepsilon)&=\left(\frac{b_1+6\varepsilon}{5+12\varepsilon}-c_2\right)\frac{5+12\varepsilon}6+2b= \frac{b_1}6+2+\varepsilon-\left(\frac56+\varepsilon\right)c_2<\frac{b_1}6+2+\varepsilon,\\
\end{align*}
\begin{align*}
f(\varepsilon)
&=\delta+\delta^2+(1-\delta)\varepsilon-\frac{40-4\sqrt{2}}{49}\left(2\delta+2\delta^2+\frac{2\delta^3}{1-\delta}+\varepsilon\right)<\\
&<\delta+\delta^2+(1-\delta)\varepsilon-\frac{34}{49}\left(2\delta+2\delta^2+\frac{2\delta^3}{1-\delta}+\varepsilon\right)=\\
&=-\frac27\delta+\frac{15}{49}\varepsilon-\frac5{49}\delta=-\frac27\delta+\frac{15}{49}\left(\varepsilon-\frac{\delta}3\right) \le-\frac27\delta;\\
v(\varepsilon)&=-c_2\left(\frac{2\delta}{1-\delta}+\varepsilon\right)+2<2.
\end{align*}
In the right-hand side of formula (\ref{Formula otsenka summi B2-2}), using condition \eqref{Formula neravenstva y}, and then substituting the obtained upper estimates for the functions $g(\varepsilon)$,  $h(\varepsilon)$, $f(\varepsilon)$ and $v(\varepsilon)$, we successively obtain
\begin{align*}
{\mathscr B}_2&\ll yx^{g(\varepsilon)}{\mathscr L}^{h(\varepsilon)} +yx^{f(\varepsilon)}{\mathscr L}^{v(\varepsilon)}\ll yx^{-\frac1{10000}}{\mathscr L}^{\frac{b_1}6+2+\varepsilon} +yx^{-\frac27\delta}{\mathscr L}^2\ll\\
&\ll y{\mathscr L}^{-A}+y{\mathscr L}^2\exp\left(-\frac27\delta{\mathscr L}\right).
\end{align*}
Using condition (\ref{Formula neravenstva y}), we obtain
\begin{align*}
T_0&=\left(xy^{-1}+\lambda x^2\right)q^{\frac12}{\mathscr L}^{A+3}\le\left(\frac{x\sqrt{q}}y+\frac{x^2}{\sqrt{q}\tau}\right){\mathscr L}^{A+3}\le
\\ &\le\left(\frac{x\sqrt{q}}y+\frac{x^3{\mathscr L}^{b_1}}{y^3\sqrt{q}}\right){\mathscr L}^{A+3}\le \left(\frac{x{\mathscr L}^b}y+\frac{x^3{\mathscr L}^{b_1}}{y^3}\right){\mathscr L}^{A+3}\le\\
&\le \left(x^{\frac1{9-4\sqrt{2}}}{\mathscr L}^{b-c_2}+x^{\frac3{9-4\sqrt{2}}}{\mathscr L}^{b_1-3c_2}\right){\mathscr L}^{A+3}=\\
&=\left(x^{\frac{9+4\sqrt{2}}{49}}{\mathscr L}^{b-c_2}+x^{\frac{27-12\sqrt{2}}{49}}{\mathscr L}^{b_1-3c_2}\right){\mathscr L}^{A+3}<\\
&<\left(x^{\frac{15}{49}}{\mathscr L}^{b-c_2}+x^{\frac{45}{49}}{\mathscr L}^{b_1-3c_2}\right){\mathscr L}^{A+3}<
0.01x.
\end{align*}
With the help of this estimate we estimate the parameter $\delta=\delta(q,T_0)$ from below. We have
\begin{align*}
\delta(q,T_0)
 &\ge\frac{c_1}{\max\left(b\ln{\mathscr L},\ \left({\mathscr L}\ln{\mathscr L}\right)^\frac{3}{4}\right)}\ge c_1{\mathscr L}^{-0.76}.
\end{align*}
Therefore
\begin{align*}
{\mathscr B}_2&\ll y{\mathscr L}^{-A}+y{\mathscr L}^2\exp\left(-\frac27\delta{\mathscr L}\right)\ll y{\mathscr L}^{-A}.
\end{align*}
\subsection{Nontrivial estimate of $S_2(\alpha;x,y)$ in the case $\lambda\le\left(8\pi y^2\right)^{-1}$}
From the found estimates for ${\mathscr B}_1$ and ${\mathscr B}_2$, taking into account formula (\ref{Formula otsenka summi B1}), we have
$$
|W(\alpha;x,y)|\ll y{\mathscr L}^{-A}.
$$
Substituting the last estimate into the right-hand side of formula (\ref{formula Sk(alpha;x,y)=int+W1+W+R2-3}), we find
\begin{align*}
S_2(\alpha;x,y)=\frac1{\varphi(q)}\sum_{\substack{n=1\\(a,q)=1}}^qe\left(\frac{an^2}q\right)\int_{x-y}^xe(\lambda u^2)du+ O\left(\frac{y}{{\mathscr L}^A}\right),
 \end{align*}
Passing to estimates and using the trivial estimate of the sum over $n$, that is, the number of terms, we find
\begin{align*}
S_2(\alpha;x,y)\ll \left|\int_{x-y}^xe(\lambda u^2)du\right|+\frac{y}{{\mathscr L}^A}.
 \end{align*}
Estimating the trigonometric integral by the magnitude of the first derivative, and using the condition $\lambda>\left(8\pi y^2\right)^{-1}$, and then the relation (\ref{Formula neravenstva y}), that is, the inequality $ y\le x{\mathscr L}^{-A-0.5b-10}$, we have
\begin{align*}
S_2(\alpha;x,y)\ll\frac1{\min\limits_{x-y\le u\le x}{2\lambda x}}+\frac{y}{{\mathscr L}^A}\ll y\left(\frac{y}x+\frac1{{\mathscr L}^A}\right)\ll \frac{y}{{\mathscr L}^A}.
 \end{align*}

\end{document}